\documentclass[12pt]{article}

\usepackage{amsmath,amsthm,amsfonts,latexsym,amscd}

\usepackage[mathscr]{eucal}

\usepackage{graphicx}

\theoremstyle{plain}
          \newtheorem{theorem}{Theorem}[section]
          \newtheorem{lemma}[theorem]{Lemma}
          \newtheorem{corollary}[theorem]{Corollary}

\theoremstyle{definition}
          
          \newtheorem{example}[theorem]{Example}

          \newtheorem{remarks}[theorem]{Remarks}

\numberwithin{equation}{section}

\addtocounter{section}{-1}

\newcommand{\gs}{\sigma}
\newcommand{\ga}{\alpha}
\newcommand{\gb}{\beta}
\newcommand{\gl}{\lambda}
\newcommand{\gm}{\gamma}

\newcommand{\ta}{\theta}

\newcommand{\xv}{\mathbf x}
\newcommand{\yv}{{\mathbf y}}
\newcommand{\tv}{\mathbf t}

\newcommand{\ptspp}[1]{p_{s_{#1},\tv_{#1}}^+}

\newcommand{\bt}[1]{b_{\tv_{#1}}}
\newcommand{\ptsp}{p_{s,\tv}^+}

\newcommand{\ptsm}{p_{s,\tv}^-}

\newcommand{\ptspm}{p_{s,\tv}^{\pm}}

\newcommand{\rn}[1]{$\mathbb R^{#1}$}

\newcommand{\conv}{\operatorname{conv}}

\newcommand{\lan}{\langle}
\newcommand{\ran}{\rangle}

\newcommand{\spk}[1]{(s_{#1},\tv_{#1})}
\newcommand{\gat}{r_\ta}

\DeclareMathOperator{\dis}{Disc}

\begin{document}
\title{The Spectral Scale and the $k$--Numerical Range}
\author{Charles A. Akemann and Joel Anderson
\thanks{The second author was partially
supported by the National Science Foundation during the period of
research that resulted in this paper.}}

%
%
%

\maketitle
\begin{abstract}
Suppose that $c$ is a linear operator acting on an $n$-dimensional complex
Hilbert Space $H$, and let $\tau$ denote the normalized trace on $B(H)$. 
Set $b_1 = (c+c^*)/2$ and $b_2 = (c-c^*)/2i$, and write $B$ for the the
spectral scale of $\{b_1, b_2\}$ with respect to $\tau$.  We show that
$B$ contains full information about
$W_k(c)$, the $k$-numerical range of $c$  for each $k = 1,\dots,n$.  We
then use our previous work on spectral scales to prove several new facts
about $W_k(c)$.  For example, we show in Theorem 3.4 that the point $\gl$
is a singular point on the boundary of $W_k(c)$ if and only if $\gl$ is an
isolated extreme point of $W_k(c)$. In this case $\gl = (n/k)\tau(cz)$,
where $z$ is a central projection in in the algebra generated by $b_1$,
$b_2$ and the identity. We show in Theorem 3.5, that $c$ is normal if and
only if $W_k(c)$ is a polygon for each $k$.  Finally, it is shown in
Theorem 5.4 that the boundary of the $k$-numerical range is the finite
union of line segments and curved real analytic arcs.
\end{abstract}

\section{Introduction and Notation}
      \medskip
The spectral scale was introduced by the present authors and Nik
Weaver in \cite{AAW} and further developed by the authors in \cite{Geom
II}. It is defined for any finite set of self-adjoint operators in a
finite von Neumann algebra. The main theme in \cite{AAW} and \cite{Geom
II} is that full spectral information about real linear combinations of
such operators is contained in the spectral scale and that much of this
information is reflected by the geometry of the spectral scale.

In the present paper we view the spectral scale from another perspective
and thereby show that full information about the $k$--numerical range of a
finite dimensional operator is also contained in its spectral scale. Thus,
in this restricted case at least, two rather different sets of data are
combined into one three dimensional, compact, convex set. Although we
restrict consideration in this paper to the finite dimensional world, we
shall write this paper using the language of operators (rather than
matrices) to emphasize how easily many of the concepts generalize to
infinite dimensional situations.

The notation developed here will be used throughout the rest of the paper.
Let $c$ denote a linear operator on an $n$--dimensional complex Hilbert
Space $H$ and let $\tau$ denote the normalized trace on $B(H)$, the
algebra of all linear operators on $H$. Write $N$ for the subalgebra 
of $B(H)$ generated
by $c, c^*$ and the identity $1$ of $B(H)$ and set $b_1 = (c+c^*)/2$ and $b_2
= (c-c^*)/2i$.   The spectral scale was defined in \cite{AAW} and \cite{Geom
II} via the map $\Psi$ defined by the formula
\[
\Psi(a) = (\tau(a),\tau(b_1a),\tau(b_2a)),
\]
and the spectral scale $B = B(b_1,b_2)$ was defined as
\[
B = \{\Psi(a):  a\in B(H), 0 \le a \le 1\}.
\]
It is convenient for the present paper to view the second and third 
real coordinates as a
single complex number.  The definition of $\Psi$ then becomes
\[
\Psi(a) = (\tau(a),\tau(ca))
\]
and we now define the spectral scale  $B = B(c)$ of $c$ by the formula
\[
B =\{\Psi(a) : a \in B(H), 0 \le a \le 1\}.
\]
Thus, we now view $B$ as a subset of $\mathbb R\times \mathbb C$, rather
than as a subset of \rn{3}. As shown in \cite [Theorem 2.4] {AAW}, $B
=\Psi(N_1^+)$, where  for any self-adjoint subalgebra $M$ of $B(H)$,
\[
M_1^+ = \{a\in M: 0\le a\le 1\}.
\]

As noted above, the geometry of the spectral scale reflects spectral data
for real linear combinations of $b_1, b_2$, i.e. matrix pencil information.
Information of this sort is widely valued as documented in \cite{Knud}
The basis for this current work is the observation that $B$ also
essentially contains the $k$-numerical range of $c$ for each $1 \le k \le
n$.

Since we shall show that several of the currently known properties of the
numerical range generalize to the $k$--numerical range, let us begin 
by reviewing these
properties. Recall that the {\bf numerical
range} of $c$ is by definition
\[
W(c) = \{\lan c\xv,\xv\ran: \xv\in H \text{ and } \|\xv\| =1\}.
\]
We use $\gs(c)$ to denote the spectrum of $c$.

\begin{theorem} The following statements hold.
\begin{enumerate}

\item $W(c)$ is a compact, convex subset of $\mathbb C$.

\item If $a= a^*$ and
\[
\ga^- = min\{\ga\in \gs(a)\} \text{ and } \ga^+ = max\{\ga\in
\gs(c)\},
\]
then $ W(c) =[\ga^-,\ga^+]$.

\item If $c$ is a normal operator, then $W(c)$ is the convex hull of its
eigenvalues.

\item  $W(c)$ is a line segment if and only if $c = \gl a + \mu 1$, where
$\gl$ and $\mu$ are complex numbers and $a$ is self-adjoint.

\item The boundary of $W(c)$ is the  union of
a finite number of  analytic arcs.

\item If $\gl\in \gs(c)$ and $\gl$ lies in the relative boundary of $W(c)$,
then $\gl$ is a reducing eigenvalue of $c$.

\item The boundary of $W(c)$ is contained in the real zero set of an
algebraic curve.

\end{enumerate}

\end{theorem}
\begin{proof}  Assertion $(1)$ is known as the {\em Toeplitz--Hausdorff
Theorem} (see \cite{Toep} and \cite{Haus}). A proof of this and the next 3
assertions can be found in \cite[\S 1.2]{Horn}. The assertion in (5) was
proved independently by  Agler \cite[Theorem 4.1]{Ag} and Narcowich
\cite[Corollary 3.5]{Nar}. Finally assertions (6) and (7) are due to
Kippenhahn \cite{Kip}.
\end{proof}

Observe that if $W(c)$ is a line segment, then it follows from  parts (2)
and (4) of Theorem 0.1 that the endpoints of $W(c)$ are reducing eigenvalues
of $c$.  Also, $W(c) = \{\gl\}$ if and only if $c = \gl 1$. In fact if $W(c)$
is a line segment and  $\xv$ is a vector in $H$ such that $\lan c \xv,\xv\ran
= \eta$ lies on an end point of $W(c)$, then  $c\xv = \eta \xv$ and $c^*\xv =
\bar \eta \xv$.  To see this, note that by part (4) of Theorem 0.1, it
suffices to show this when $W(c)$ is contained in $[0,\infty)$  and $\eta
=0$.  In this case we get that $c \ge 0$ and so it has a positive square root.
Hence, we get
\[
\|\sqrt{c}\xv\|^2 = \lan \sqrt{c}\xv,\sqrt{c}\xv\ran = \lan c\xv,\xv\ran = 0
\]
and so $c\xv = (\sqrt{c})^2\xv = 0$.

The numerical range  has a {\bf corner} at $\gl$  if $W(c)$ has dimension two
and there is more than one tangent line of support for $W(c)$ at $\gl$.   The
presence of a corner in the numerical  range signals the fact that $c$ enjoys
a special structure.  Further, it turns out that  corners in $W(c)$ must be
{\bf lineal} in the sense that there are two tangent lines of  support which
intersect the boundary in line segments of positive length.  Since we  shall
show that an analogous results hold for the $k$-numerical range, we now
discuss  this in more detail. Recall that a point on the boundary of a convex
subset of
\rn{2} is said to be {\bf singular} if the the boundary curve is not
differentiable at this point.

\begin{theorem} If $W(c)$ has dimension two and the boundary of $W(c)$ is
singular at $\gl$, then the following statements hold.
\begin{enumerate}
\item $W(c)$ has a corner at $\gl$.
\item $\gl$ is a reducing eigenvalue for $c$.
\item  $\gl$ is lineal.
\end{enumerate}
\end{theorem}

\begin{proof} Let us begin by presenting a proof of the assertion in 
(1). Rotating and
translating if necessary, we may assume that  $\gl = 0$,
$W(c)$ lies in the upper half plane and that the positive imaginary
axis intersects
the interior of $W(c)$. In this case we may find a convex function
$f$ defined on an open interval containing 0 whose graph gives a portion
of the boundary of $W(c)$ which contains 0.  Since 0 is a singular point of
the boundary, $f$ is not differentiable at 0 and since $f$ is convex, it
follows that
\[
f_{-}^\prime(0) = \lim_{h\to 0^-}\frac{f(h))}{h} < \lim_{h\to
0^+}\frac{f(h))}{h} =
f_{+}^\prime(0).
\]
Thus, the lines though $0$  with slopes $f_{\pm}^\prime(0)$ are
tangent to $W(c)$.

Kippenhahn established assertion (2) in \cite[Satz 13]{Kip}.  The 
proof offered below
seems to be new.  It is convenient to rotate once more so that $\gl = 
0$ and $W(c)$ lies
in the right half plane so that the corner has the form shown below 
where both $L_i$ are
tangent to $W(c)$ at 0, $L_1$ has positive slope and $L_2$ has negative slope
as shown below.

\vskip .2in
\setlength{\unitlength}{1 in}
\begin{picture}(2,2)(-2.5,0)
\put(-1.5,0){\line(1,0){2}}
\put(-1.3,-1.25){\line(0,1){2.5}}
\put(-1.3,0){\line(1,1){1}}
\put(-1.3,0){\line(1,-1){1}}
\put(.55,-.03){$x$}
\put(-1.33,1.32){$y$}
\put(-.28,1){$L_1$}
\put(-.28,-1.1){$L_2$}
\put(-1.34,-.04){$\bullet$}
\put(-1.75,.05){$\gl=0$}
\end{picture}
\vskip 1.8in
Since $W(c)$ is contained in the right half plane, we get $b_1 \ge 0$
and since there is a unit vector $\xv$ such that $\lan b_1\xv,\xv\ran =
0$, we get  $b_1\xv = 0$ by the remark following Theorem 0.1.

We may now select $\ta \ne 0$ so small that $W(e^{i\ta}c)$ is also
contained in the right half plane.  Arguing as above, we get that
\[
\text{Re}(e^{i\ta}c)\xv = (\cos\ta b_1 - \sin\ta b_2)\xv =  0
\]
and so $b_2\xv = 0$ and therefore $b_1\xv = b_2\xv = c\xv = 0$.
Hence, $\gl$ must be a reducing eigenvalue for $c$.

Assertion $(3)$ is due to Lancaster \cite[Corollary 4]{Lanc}.
\end{proof}
Observe that the proof above shows that if $C$ is any compact, convex, two
dimensional subset of \rn{2}, then the boundary of of $C$ is singular at $\gl$
if and only if $C$ has a corner at $\gl$.

The $k$--numerical range of $c$ is defined by the formula
\[
W_k(c) = \left\{\frac{1}{k}\sum_{i=1}^k \lan
c\xv_i,\xv_i\ran:\text{ the $\xv_i$'s are orthonormal}\right\},\quad 1
\le k \le n.
\]
Observe that when $k = n$ we have
\[
W_n(c) = \tau(c)
\]
so that $W_n(c)$ consists of a single point.  As the results in this 
paper make clear
it is natural to also include the case where $k = 0$ and to define 
$W_0(c) = 0$.
Since we have $W_1(c) = W(c)$,  this notion is a generalization of the
standard numerical range.  Let us now review the basic  properties of
the $k$-numerical range.

\begin{theorem}
If $1 \le k \le n$, then the following statements hold.
\begin{enumerate}

\item $W_k(c)$ is a compact, convex  subset of
$\mathbb C$.

\item  $\displaystyle W_k(c) =
\left\{\frac{n}{k}\tau(cp):\text{
$p$ is a projection of rank $k$}\right\}$.

\item We have $kW_k(c) = n\tau(c) - (n-k)W_{n-k}(c)$.

\item If $a$ is self-adjoint with eigenvalues $\ga_1 \ge \cdots
\ge\ga_n$, $\ga_k^+ =\ga_1 + \cdots + \ga_k$ and $\ga_k^- = \ga_{n-k+1} +
\cdots + \ga_n$ then $\displaystyle W_k(a) = [\ga_k^-,\ga_k^+]$

\item If $\gb_k^+$ denotes the sum of the $k$ largest eigenvalues
of
$b_1$, then the line $x = \gb_k^+/k $  is tangent to $W_k(c)$.

\end{enumerate}
\end{theorem}
\begin{proof} The assertion in $(1)$ is due to Berger who introduced 
the $k$--numerical
range in his thesis \cite{Berg}. A proof may be found in 
\cite[Problem 167]{Halmos}.
The assertion in $(2)$ follows from part $(1)$ and a simple 
calculation.  Assertion $(3)$
follows from the fact that $\tau(cp) + \tau(c(1-p)) = \tau(c)$.  The proof of
assertion $(4)$ is straight forward.  For assertion $(5)$ observe that
\[
\frac{1}{k}\sum_{i=1}^k \lan(c\xv_i,\xv_i\ran = \frac{1}{k}\sum_{i=1}^k
\lan(b_1\xv_i,\xv_i\ran +i\frac{1}{k}\sum_{i=1}^k
\lan(b_2\xv_i,\xv_i\ran
\]
so that $\gb_k^+ + i\gm \in W_k(c)$ for some real $\gm$.  Further, since
$\gb_k^+$ is the sum of the $k$ largest eigenvalues of $b_1$, if
$\gb+i\delta\ \in W_k(c)$, then  $\gb \le \gb_k^+$.
\end{proof}

Let us now describe our results in more detail.  The key to
understanding the role played by the $k$-numerical range in the
spectral scale is the notion of an {\em isotrace slice} of the spectral
scale.  If $0 \le t \le 1$, then the {\bf isotrace slice}  of $B$ at $t$ is
by definition
\[
I_t = \{\xv = (x_0,z)\in B: x_0 = t\}.
\]
We prove in Theorem 1.3 that if $0 < k <  n$ and we define the map
$\pi_k$  from  $\mathbb C$ to $\mathbb R\times \mathbb C$ by  $\pi_k(z) =
(k/n,kz/n)$,  then $\pi_k$ is an affine map that is a bijection from
$W_k(c)$ onto $I_{k/n}$.  Thus, we may view the $k$-numerical range as a
subset of the spectral scale.

It is also shown in section 1 (Theorem 1.1) that the extreme points of $B$
lie on the isotrace slices of the form $I_{k/n}$ for $k = 0,1,\dots,n$ and
so $B$ is the convex hull of this finite collection of sets.  In section 2 we
present some examples (and pictures) of various spectral scales.  Section
3 contains our results on corners in the $k$-numerical range.  It is shown
that several of the known facts about corners on the boundary of $W(c)$
generalize to $W_k(c)$ by using the additional structure provided by the
spectral scale.  For example, in Theorem 3.4, we show that a singularity on
the boundary of $W_k(c)$ always occurs at an isolated extreme point of
$W_k(c)$ and  such points correspond to central projections in the algebra
$N$, i.e. reducing subspaces for $c$. We also show that $c$ is normal if and
only if $W_k(c)$ is a polygon for $0<k< n$ (Theorem 3.5).

Sections 4 and 5 are devoted to establishing that the boundary of
the $k$-numerical range is the finite union of line segments and curved
real analytic arcs (Theorem 5.4). In section 4 we review some classical
background which is required for our analysis and then use this to derive
our results in section 5. Section 6 contains some open question, stated as
conjectures.

As our results show, the spectral scale provides a new way to study
$n$-tuples of self-adjoint finite dimensional operators.  In the case
under study here, when there are just two operators so that $B$ is a
subset of three dimensional real euclidean space (or $\mathbb R\times
\mathbb C$), we may actually visualize $B$ as shown in the examples and
pictures in section 2.  These pictures were created using a MATLB program
written by Jeff Duzak as part of an REU research project supervised by the
second author. Readers can contact  the second author for a copy of this
program which is  quite useful for testing conjectures.

\section{Isotraces and extreme points}
\bigskip
  The spectral scale has a striking structure in the finite 
dimensional case under
consideration here which we describe in the next proposition.

\medskip
\begin{theorem} If $\xv$ is an extreme point of $B$, then it lies in 
an isotrace slice of
the form $I_{k/n}$ where $k = 0,1,\dots,n$.
\end{theorem}
\begin{proof} Since $\xv$ is an extreme point of $B$, it has the form 
$\Psi(p)$ where
$p$ is a projection in $N$ by \cite[Theorem 2.3(1)]{AAW}.  If $p$ has 
rank $k$, then
$\tau(p) = k/n$ and so $\Psi(p) = (k/n,s)$, where $s = \tau(cp)$. 
Hence, $\xv \in
I_{k/n}$.
\end{proof}

For example, if $n=3$,  then extreme
points of $B$ come from projections of trace
$0$, $1/3$, $2/3$ or $1$. While $I_0$ and $I_1$ are always the single
points $0$ and $\Psi(1)$, in generic examples  $I_{1/3}$ and 
$I_{2/3}$ are solid
ellipses whose boundaries consist of extreme points of $B$. Hence, 
generically, the
boundary of $B$ between two successive isotrace slices consists of one
dimensional faces.  The boundary between $I_0$ and $I_{1/3}$  is
typically a skewed cone, and, by the symmetry of $B$, the same is true for
the boundary of $B$ between $I_{2/3}$ and $I_1$. However, as one can
see in Examples 2.2 and 2.3 below, $B$ may have planar faces.  A complete
description and interpretation of the faces of  $B$ may be found in
\cite[\S 3]{Geom II}. Further,  in \cite[Corollary 5.4]{Geom II}
we showed that $N$  is abelian and finite-dimensional  if and only if
the spectral scale has a finite number of extreme points.  Thus, in
finite dimensions, one can ``see'' that $N$ is abelian from the shape
of $B$.

We now show how the $k$-numerical range may be identified with the
isotrace slice $I_{k/n}$.  This identification depends on a simple
convexity result, which we now present.  If $0 < t < 1$ and $M$ is
any self-adjoint
subalgebra of $B(H)$, then write
\[
M_{1,t}^+ = \{a\in M_1^+: \tau(a) = t\}.
\]

\begin{lemma} The extreme points of $B(H)_{1,k/n}^+$ are precisely the
projections of rank $k$.
\end{lemma}
\begin{proof} If $p$ is a projection of rank $k$, then it  an
extreme point of $B(H)_1^+$ and so it is also an extreme point of
$B(H)_{1,k/n}^+$. For the converse suppose $a$ is in $N_{1,k/n}^+$, 
but $a$ is not a
projection and write  $\ga_1,\dots,\ga_n$ for the  eigenvalues of 
$a$. Since  $\tau(a) =
k/n$ we get that $\ga_1+\cdots + \ga_n = k$ and  since $a$ is not a 
projection  we have
$0 <
\ga_i < 1$ for at least one index $i$.  Since $\ga_1+\cdots + \ga_n = 
k$, there must also
be an index $j\ne i$ such that $0 < \ga_j < 1$.  Relabeling if necessary, we 
may
assume that $0 <\ga_1 \le \ga_2 < 1$.

       Since $\ga_1+\ga_2 -1<\ga_1$ and $0 < \ga_1$ we may  select $\gl$ and
$\gm$  such that
\[
\max\{0,\ga_1+\ga_2 - 1\} < \gl < \ga_1 \text{ and } \ga_2 < \gm <
\min\{1,\ga_1+\ga_2\}.
\]
Now write
\[
a_1  =
\bmatrix
\gl&0&0&\hdots&0 \\
0&\ga_2+ \ga_1 -\gl&0&\hdots&0\\
0&0&\ga_3&\hdots&0 \\
\vdots&\vdots&\vdots&\ddots&\vdots\\
0&0&0&\hdots&\ga_n
\endbmatrix
\text{ and }
a_2 =
\bmatrix
\gm&0&0&\hdots&0 \\
0&\ga_2+ \ga_1 -\gm&0&\hdots&0\\
0&0&\ga_3&\hdots&0 \\
\vdots&\vdots&\vdots&\ddots&\vdots\\
0&0&0&\hdots&\ga_n
\endbmatrix.
\]
Observe that $a_1 \ne a_2$ because $\gl < \ga_1 \le \ga_2 < \gm$.  Since
the diagonal entries of $a_1$ and $a_2$ lie in $[0,1]$  and sum to
$k$, these  matrices are elements of  $B(H)_{1,k/n}^+$.  As $\gl < \ga_1 <
\gm$ there is a real number $t$ with $0 < t < 1$, such that  $\ga_1 =
t\gl + (1-t)\gm$.  Next, note that
\begin{align*}
t(\ga_1 + \ga_2 -\gl) + (1-t)(\ga_1+\ga_2 +\gm) &= \ga_1 + \ga_2 -(t\gl
+(1-t)\gm)\\ &= \ga_1+ \ga_2 - \ga_1 = \ga_2
\end{align*}
and therefore  $a = ta_1 + (1-t)a_2 $.  Thus $a$ is not an extreme
point.
\end{proof}

\begin{theorem} If $0 < k < n$ and we define the map $\pi_k$ of $\mathbb
C$ into $\mathbb R\times \mathbb C$ by the formula
\[
\pi_k(\gl) = (k/n,\gl),
\]
then $\pi_k$ is an affine  bijection from
$W_k(c)$ onto $I_{k/n}$.
\end{theorem}
\begin{proof}
If $\gl \in W_k(c)$, then by part$(2)$ of Theorem 0.3  there is a
projection $p$ with rank $k$ such that
\[
\gl = \frac{n}{k}\tau(pc)
\]
and the point $(\tau(p),\tau(pc))$ is in $I_{k/n}$ because $\tau(p) =
k/n$. Thus,
\[
\pi_k(\gl) = (k/n,k\gl/n) = (k/n,\tau(pc)) \in I_{k/n}.
\]
So, $\pi_k$ maps $W_k(c)$ into $I_{k/n}$, and $\pi_k$ is clearly a
one-to-one map.

Now suppose that $(k/n,\gl) \in I_{k/n}$ so that $\gl = \tau(ca)$ for some
$a\in N_{1,k/n}^+$. By Lemma 1.2, the Krein--Milman Theorem  \cite[ Theorem
2.6.16 ]{Web} and \cite[Theorem 17.1]{Rock} we have that $a$ is a convex
combination of
projections of rank $k$. Since $W_k(c)$ is convex by part $(1)$ of Theorem
0.3 and points of the form $(n/k)\tau(cp)$ are in $W_k(c)$ by part $(2)$,
we get that $(n/k)\gl \in W_k(c)$.
\end{proof}

\bigskip
\section{Examples}

In this section we describe the spectral scale for four examples.

\example

If
\[
c =
\bmatrix
1+\dfrac{i}{2}&\dfrac{i}{2}\\[5pt]
\dfrac{i}{2}&\dfrac{i}{2}
\endbmatrix,
\]
then the nontrivial projections in $N$ must have trace $1/2$ and so
$I_{1/2}$ is the only isotrace slice of $B$ that contains nontrivial
extreme points.  It turns out that $I_{1/2}$ is the disk of radius
$1/4$  centered at $(1/2,(1+i)/4)$.   Since
\[
B = \conv(0,I_{1/2},(1,(1+i)/2))
\]
the spectral scale in this case is a pair of skewed circular cones
joined at their bases as shown below.

\begin{figure}[h]
\centering
\includegraphics[scale=.75]{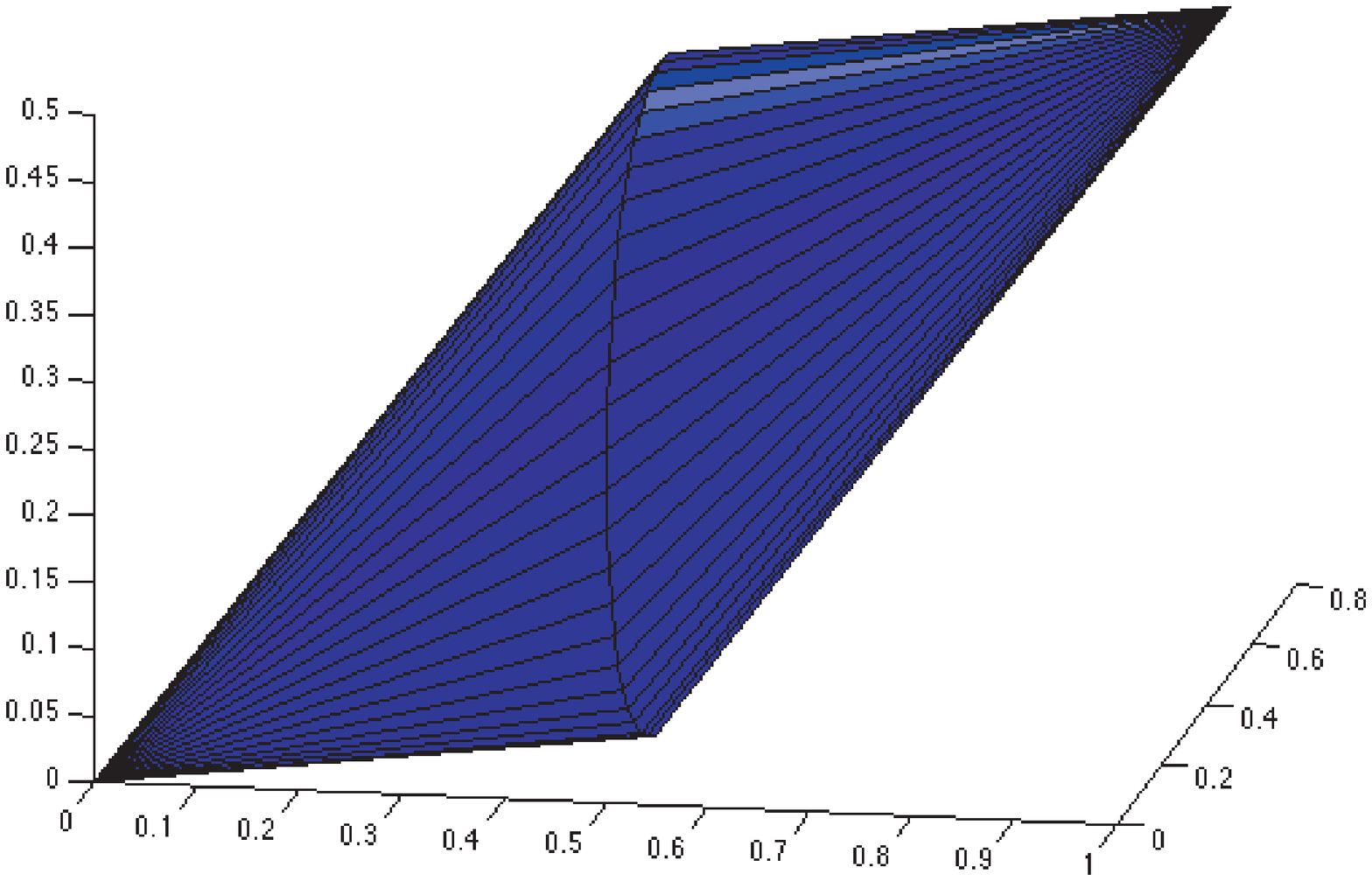}
\end{figure}

In this case the extreme points of $B$ are $0,\Psi(1)$ and points on
the circular boundary of $I_{1/2}$.  The one dimensional
faces are the line segments joining $0$ and $\Psi(1)$ to the extreme
points on the circle. There are no faces of dimension two.

\endexample
\newpage

We now present a $3\times 3$ example where the spectral scale has a
``flat spot''.

\example  Write
\[
b_1 =
\bmatrix
1&0&1\\
0&2&1\\
1&1&3
\endbmatrix,
b_2 =
\bmatrix
1&0&0\\
0&1&0\\
0&0&0
\endbmatrix
\text{ and }
c = b_1 + ib_2 =
\bmatrix
1+i&0&1\\
0&2+i&1\\
1&1&3
\endbmatrix.
\]
We have that $\Psi([0,b_2])$ is a face of $B$ by \cite[Theorem
2.3]{AAW} (here if $0 \le a^- \le a^+$, then $[a^-,a^+] = \{a: a^-\le
a\le a^+\}$).  Further it  follows from the results in \cite[\S3]{Geom
II} that this face is two dimensional.  In fact as can be seen from the
figure below, it is diamond shaped.

\begin{figure}[h]
\centering
\includegraphics[scale=.85]{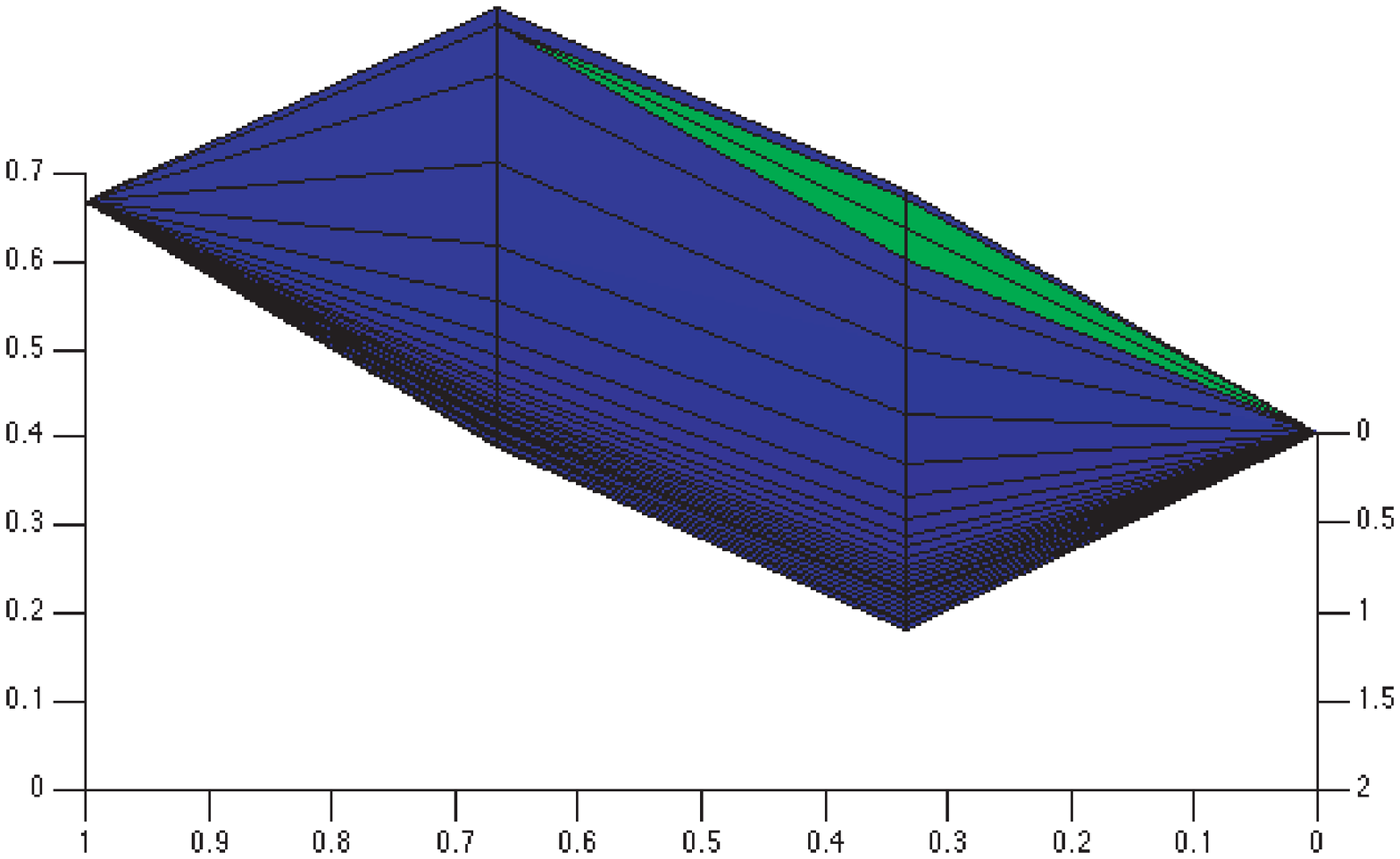}
\end{figure}

\endexample

\newpage

\example  The facial structure of the spectral scale can be quite complex.
The next example there are eight faces of dimension two.
If we write
\[
b_1 =
\bmatrix
0&0&0&0&0\\
0&0&0&1&0\\
0&0&0&0&1\\
0&1&0&0&0\\
0&0&1&0&0\\
\endbmatrix,
b_2 =
\bmatrix
0&0&0&1&0\\
0&0&0&0&1\\
0&0&0&0&0\\
1&0&0&0&0\\
0&1&0&0&0\\
\endbmatrix
\text{ and } c = b_1 + ib_2,
\]
then the spectral scale is as shown below.

\newpage
\begin{figure}[h]
\centering
\includegraphics[scale=.85]{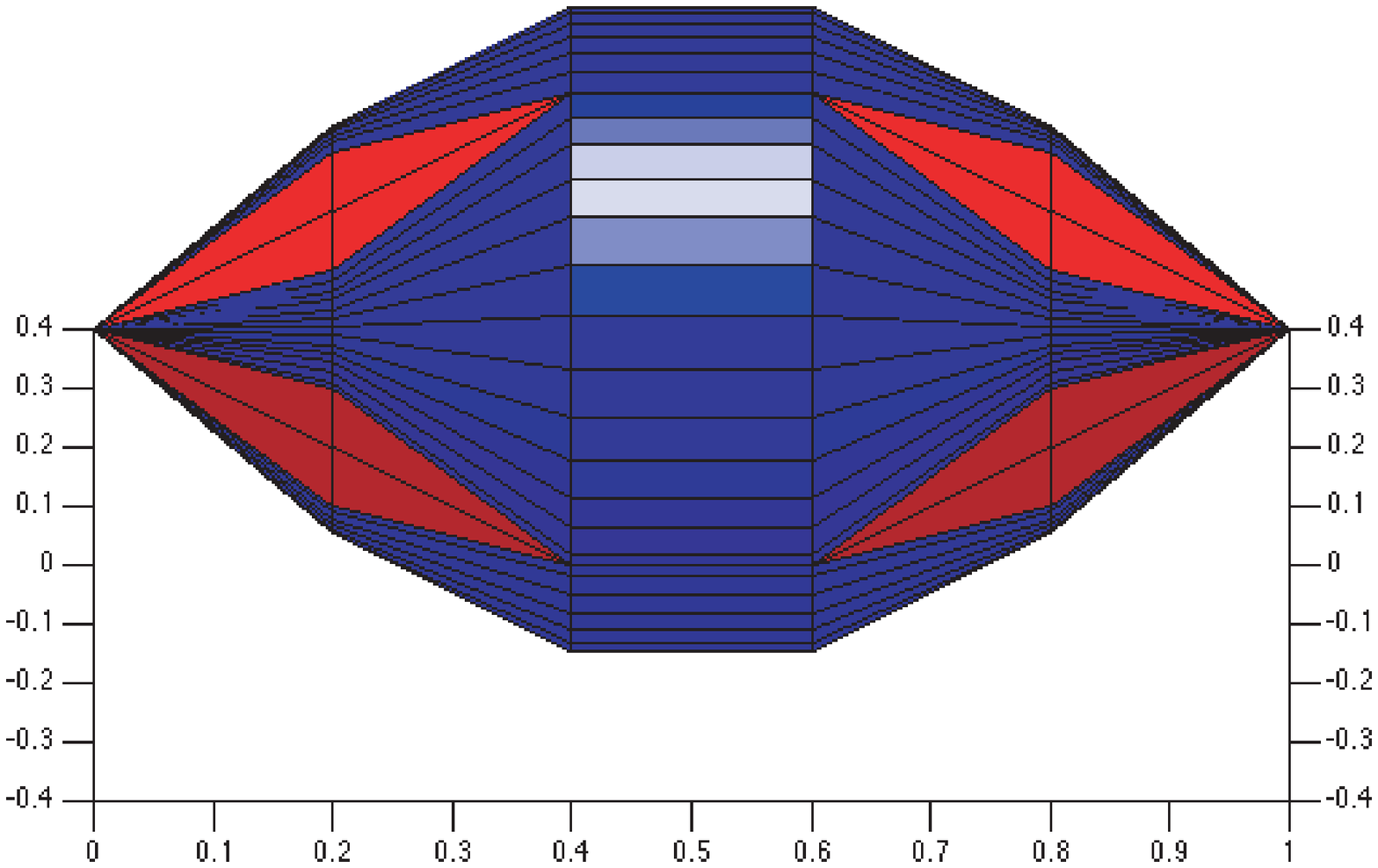}
\end{figure}

Although it appears from the figure that the lighter portions between the
isotrace slices $I_{2/5}$ and $I_{3/5}$  are  two dimensional
faces, in fact they consist of faces of dimension one. This occurred as a
result of the MATLAB shading routine employed to draw this picture.
\endexample

\example  In \cite[Example 3.5]{AAW} we showed that in the
noncommutative case the
spectral scale is not a complete invariant for the operators $b_1$ and $b_2$ by
exhibiting two pairs of inequivalent  self-adjoint $3\times 3$
matrices that shared the
same spectral scale.  In this example each pair generated the full
von Neumann algebra
$B(H)$.  In our final example we show that the same phenomenon can
occur even if one
pair is reducible.  Write
\[
c_1 =
\bmatrix
1+i&0&0\\
0&1+2i&1\\
0&1&1
\endbmatrix
\text{ and }
c_2 =
\bmatrix
1&\dfrac{1}{\sqrt2}&0\\
\dfrac{1}{\sqrt2}&1+i&\dfrac{1}{\sqrt2}\\
0&\dfrac{1}{\sqrt2}&1+2i
\endbmatrix.
\]
It is clear that for $c_2$ the algebra $N$ is the full algebra
$B(H)$, while this is not true for $c_1$ because $N$ has non-trivial
center.

On the other hand, the associated spectral scales are equal.  This
follows from the fact that each of these matrices has the same
characteristic polynomial and part $(2)$ of \cite[Theorem 3.2]{AAW}.

%
%

In fact even more is true.  Write $S(z) = \text{Re}c_1 + z \text{Im}c_1$
and $T(z) = \text{Re}c_2 + z \text{Im}c_2$.  A straight forward
calculation shows
\[
\det(S(z) -w1) = \det(T(z)-w1) = (-2zw+2z+w^2-2w)(w-1-z)
\]
so that even if $N$  is irreducible, the associated characteristic
polynomial  may be
reducible.

\endexample

\bigskip

\section{Corners in isotrace slices and the $k$-numerical range}
\bigskip

     In this section we show that by using the spectral scale we can establish
results analogous to Theorem 0.1(6) and Theorem 0.2  for the $k$--numerical
range.  Let us begin by recalling an idea that proved useful in \cite{Geom
II}. If
$\xv$ is an extreme point of a convex set $X$ in $\mathbb R^k$ which is
isolated in the set of all extreme points of $X$, then we say that
$\xv$ is an {\bf isolated extreme point} of $X$.

Recall that a face in $B$ is said to be a {\bf sharp face} if it is
contained in at least two hyperplanes of support \cite[Definition
4.1]{Geom II}.  We showed in \cite[Lemma 5.2]{Geom II} that an isolated
extreme point of $B$ is always a sharp face. Our main result on corners is
Theorem 3.4 below.  The proof of this result requires three technical lemmas
which we present next. Since the first result is just a slight 
generalization of Lemma
1.2 and uses the notation from that lemma, the proof is somewhat abbreviated.

\begin{lemma}If $a$ is an extreme point of $N^+_{1,t}$, then there is a
projection $p$ in $N$ such that either $a = p$ or
\[
a = p + \ga q,
\]
where  $q$ is projection in $N$ which is  orthogonal to  $p$, $qNq$
has dimension one and
$0 <\ga < 1$.
\end{lemma}
\begin{proof} Fix an extreme point $a$ of $N^+_{1,t}$.  If $a = p$ is a
projection, then the proof is complete.  So, assume that this is not
the case, let $\ga_1,\dots,\ga_m$ denote the distinct eigenvalues of
$a$ and let $q_1,\dots, q_k$ denote the corresponding
eigenprojections.  Since $a$ is not a projection, there is at least one
index $i$  such that $0 < \ga_i < 1$.  If there were another index $j$
with $0 < \ga_j < 1$, then  by adding (resp. subtracting) and small
amount to $\ga_i$ and subtracting (resp. adding) a small amount to
$\ga_j$ we would get two new elements of $N^+_{1t}$ whose average is $a$
so that $a$ would not not be an extreme point of this set.  Hence, $0<
\ga_i < 1$ and all other indices are 0 or 1.

If  $q_iNq_i$ did not have dimension one, then it would contain  two
nonzero orthogonal projections and, arguing as in the previous
paragraph,  we would again get that $a$ is not an extreme point.
Hence, $a$ has the indicated form.
\end{proof}

In the following lemma it is convenient to  return (briefly) to the
picture of $B$ as lying in \rn{3} and then to use some results \cite {Geom II}
which are  stated in this
framework. us, a corner in $W_k(c)$ corresponds to a corner
in  $I_{k/n}$.  If $C$ is a convex subset of \rn{2} with a corner at
$\gl$, then $C$ admits an infinite family of tangent lines of support at
$\gl$ and these lines all lie in a sector of a disk with maximal angle.  We
call the tangent lines of support that bound this sector the {\bf sectorial}
tangent lines of support.

\begin{lemma} If $k$ is an integer with $1 \le k \le n$, $I_{k/n}$ has
dimension two and the point $\xv =(k/n,r_1,r_2)$
lies at a corner on the isotrace slice $I_{k/n}$, then the following
statements hold.
\begin{enumerate}

\item There is a sharp face $F$  in $B$ of dimension at most one that
contains $\xv$.

\item $F= \Psi([z^-,z^+])$, where $z^\pm$ are central projections.

\item The points $\Psi(z^\pm)$ are isolated extreme points of $B$.

\end{enumerate}
\end{lemma}
\begin{proof} As $I_k/n$ has dimension two, there are two distinct sectorial
tangent lines of support of $I_{k/n}$ in the plane of $I_{k/n}$ at the corner
$\xv$. Denote these tangent lines  by $L_1$ and $L_2$.  Each $L_i$  meets B
only in the boundary of the isotrace slice $I_{k/n}$ and so they are each
disjoint from the interior of B.  By
\cite[Corollary 2.4.11]{Web}, there are distinct planes $P_1$ and $P_2$
such that $L_i\subset P_i$ and each $P_i$ is disjoint from the interior
of $B$. Since each plane contains $(k/n,r_1,r_2))$, they are planes
of support for $B$. Thus, if we write $F = P_1\cap P_2\cap B$, then $F$ is a
sharp face of $B$.

Since we are now regarding $B$ as a subset of \rn{3} we
may use the  results in \cite{AAW}
and \cite{Geom II}  to get that each  hyperplane of support is determined by
a  {\em spectral pair} of the form $\spk{}$, where $s$ is a real number and
$\tv = (t_1,t_2)$ is a nonzero vector in \rn{2}.  Specifically, by
\cite[Theorem 2.3]{AAW} if $P$ is a hyperplane of support for $B$, then
there is a spectral pair $\spk{}$ such that points on$P$ satisfy an equation
of the form
\[
-sx_0 + t_1x_1 + t_2x_2 = \ga.
\]
The constant $\ga$ is determined as follows.  We write $\bt{} = t_1b_1 + t_2
b_2$  and let $\ptsp$ and $\ptsm$ denote the spectral
projections of $\bt{}$ corresponding to the intervals $(-\infty,s]$ and
$(-\infty,s)$.  With this, we have
\[
\ga = \tau((\bt{}-s 1)\ptspm).
\]
Observe that $(-s,\tv)$ is a normal vector for this plane

With this, since $P_1$ and $P_2$ are distinct faces of $B$ that have nonempty
intersection, their normal vectors are linearly independent and so there exist
linearly independent spectral pairs $(s_1,\tv_1)$ and $(s_2,\tv_2)$ such that
each $P_i$ has the equation
\[
-s_ix_0 + t_{i1}x_1 + t_{i2}x_2 = \ga_i,
\]
where $\tv_i = (t_{i1},t_{i2})$. Since the point $\xv =
(k/n,r_1,r_2)$ lies in each
plane, we get
\[
-s_i(k/n) + t_{i1}r_1 + t_{i2}r_2  = \ga_i.
\]
Now fix a point $(x_{i0},x_{i1},x_{i2})$ on the tangent line $L_i$
and observe that since
$I_{k/n}$ lies in the plane $x_0 = k/n$, we get $x_{i0} = k/n$. Since
this point also lies
in the plane $P_i$, we get
\[
-s_i(k/n) + t_{i1}x_{i1} + t_{i2}x_{i2} = \ga_i = -s_i(k/n) +
t_{i1}r_1 + t_{i2}r_2
\]
and so,
\[
t_{i1}(x_{i1}-r_1) + t_{i2}(x_{i2}-r_2) = 0.
\]
If the vectors $\tv_1$ and $\tv_2$ were linearly dependent, then it
would follow that the
tangent lines $L_1$ and $L_2$ are parallel.  Since these lines
intersect at a corner
of the isotrace slice $I_{k/n}$, this is impossible  and therefore
the $\tv_i$'s must be
linearly independent.

Applying \cite[Corollary 4.7]{Geom II}, we get that $F$ has the form
$\Psi([z^-,z^+])$, where $z^\pm$ are central projections.  Next, since
$F$ is a face in $B$, the end points $\Psi(z^\pm)$ are extreme points of $B$.
	Since the projections $z^\pm$ are central, the points $\Psi(z^\pm)$ are
isolated extreme points of $B$ by \cite[Theorem 5.4(2)]{Geom II} and so
assertions $(2)$ and $(3)$  hold.

\end{proof}

Let us now return to our previous notation so that we  regard $B$ as a
subset of $\mathbb R\times \mathbb C$. We will use the following  to 
establish $(d)$
in  part $(3)$.

\medskip

\begin{lemma}  If $G = \Psi([z^-,z^+])$ is any face of $B$ that
satisfies statements $(1),(2)$ and $(3)$ of Lemma 3.2 and $(t,\gamma) 
\in G$, then
$(t,\gamma)$ is an isolated extreme point of $I_t$.
\end{lemma}

\begin{proof}  Since $G$ is a face, its
intersection with $I_t$ is a face of $I_t$, which consists  of the
single extreme point $(t,\gamma)$ , by property (1) and the
fact that faces of $B$ are transverse to isotrace slices by \cite[Theorem
6.4(1)]{Geom II}. If $(t,z)$ were not an isolated extreme point of
$I_t$, then there would be a sequence $(t,\gamma_j)$ of extreme points
in $I_t$ that converges to $(t,\gamma)$. Since the inverse images of
extreme points contain extreme points, each of these points would have
the form $\Psi(p_j)$ or $\Psi(p_j +\ga_jq_j)$, where $p_j,q_j$ and
$\ga_j$ are as described in Lemma 3.1. We shall assume that, after
passing to a subsequence, $(t,\gamma_j)=\Psi(p_j +\ga_jq_j)$ and
further that $G$ is one dimensional. The other cases are handled
similarly, but somewhat more easily.

Since the dimension of $N$ is finite, the range of the trace on the
projections in $N$ is finite and since $\tau(p_j + \ga_j q) = t$ for
each $j$, the $\ga_j$'s also form a finite set.  So, by passing to
another subsequence, we may assume that all $\ga_j = \ga$ and $0< \ga
<1$.  We can now find another subsequence (which we continue to index
with $j$'s) such that $p_j$ converges to $p$ and  $q_j$ converges to
$q$  It follows readily that $p$ and $q$ are orthogonal projections,
and $qNq$ has dimension one.

Since $\Psi(p+\ga q) = (t,\gm)\in G$ we have
\[
z^- \le p+\ga q \le z^+
\]
and since $G= \Psi([z^-, z^+])$ is one dimensional, we  get that  $z =
z^+-z^- \ne 0$. Hence, $zb_iz$ is a scalar  multiple of $z$ for $i=1,2$
by \cite[Lemma 3.3(1)]{Geom II}. Since $z$ is central,  $z(p+\ga q)$ is
a multiple of $z$. Thus, $p+\ga q = z^- +z(p+\ga q)$ and so this
element is central.  But, since $0 <\ga < 1$   this means that both $p$
and $q$ are central.  This is impossible since the central projections
are isolated in the set of all projections.  Hence, $(t,
\gamma)$ must be an isolated extreme point of $I_t$.
\end{proof}

With these preparations, we may now present the main Theorem of this
section.    Recall that the map
$\pi_k$  from
$\mathbb C$ to
$\mathbb R\times \mathbb C$ defined by
\[
\pi_k(z) = (k/n,kz/n)
\]
is an affine map that is a bijection from $W_k(c)$ onto $I_{k/n}$ by
Theorem 1.3.

\begin{theorem}If $W_k(c)$ has dimension two and if $0 < k < n$, then the
following statements are equivalent.
\begin{enumerate}
\item $\gl$ is a singular point on the boundary of $W_k(c)$.

\item $W_k(c)$ has a {\em corner} at $\gl$.

\item There is a face $F$ of $B$ that contains the point $(k/n,k\gl/n)$
which enjoys the properties described below.
\begin{enumerate}

\item $F$ is a sharp face in $B$ of dimension at most one.

\item $F= \Psi([z^-,z^+])$, where $z^\pm$ are central projections.

\item The points $\Psi(z^\pm)$ are isolated extreme points of $B$.

\item There are faces $F_1$ and $F_2$  of dimension two such that
$F_1\cap F_2 = F$ and if $(t,\gm) \in F$, then  $F_i\cap I_{t}$ is a
line segment for each $i$. In particular, if $L_1$ and $L_2$ denote the
sectorial tangent lines of support to $W_k(c)$ at $\gl$,  then their images
$\pi_k(L_i)$ in $B$ intersect the relative interiors of the $F_i$'s.

\end{enumerate}

\item $\gl$ is an isolated extreme point of $W_k(c)$.
\end{enumerate}
\end{theorem}

\begin{proof} As in the proof of $(1)$ in Theorem 0.2 and the remark
following this theorem, we have that
$(1) \implies (2)$ because the boundary of $W_k(c)$ is convex.

Now suppose that $(2)$ holds. In this case  $\pi_k(\gl) = (k/n,k\gl/n)$ lies
at a corner of the isotrace slice $I_{k/n}$ and so we may apply  Lemma 3.2
to get that  parts $(a),(b)$ and $(c)$ of part $(3)$ are true.

Hence, if $(a), (b)$ and $(c)$ hold, then for each $(t,\gm) \in F$,
there are faces $F_{t,1}$ and $F_{t,2}$ in $B$ of dimension two such
that $(t,\gm) \in F_{t,i}$ by  \cite[Lemma
6.5(3)]{Geom II} and Lemma 3.3. If
$F$ has dimension zero so that it is a single point, then we must have
$F_{t,i} = F_i$ and so $(d)$ is true in this case.

Now suppose that $F$ has dimension one, fix $(t,\gm)$ in the relative
interior of $F$ and consider the face $F_t = F_{t,1}\cap F$.  If $F_t$
were a single point, then it would be an extreme point of $B$, which is
impossible because it lies in the relative interior of $F$.  Hence $F_t$
has dimension one.  Since $F_t$  is contained in $F$ we must have $F_t = F$.
Since it is obvious (and straightforward to prove) that three distinct faces
of dimension two in $B$ cannot intersect in a face of dimension one, we must
have that $F_{t,1}$ is either $F_1$ or $F_2$.  Hence assertion (d) holds in
this case.

Thus, in all cases we get that the faces $F_i$ have dimension two.
Since such faces are transverse to the isotrace slices of $B$ by
\cite[Lemma 6.4(1)]{Geom II}, each line $\pi_k(L_i)$ intersects $F_i$ in
its relative interior. Thus,
$(2)\implies (3)$.

Next, if $(3)$ is true, then it is clear that $\gl$ is an isolated
extreme point of $W_k(c)$ and so $(3)\implies (4)$.  If $(4)$ holds, then
it is also clear that $(1)$ holds and so these four conditions are
equivalent.
\end{proof}

\begin{theorem} The operator $c$ is normal if and only if $W_k(c)$ is a
polygon for $0 <  k < n$.
\end{theorem}
\begin{proof}Since each $W_k(c)$ is a multiple of the isotrace sliced
$I_{k/n}$ of $B$, the theorem is equivalent to showing that $B$ has a
finite number of extreme points if and only if $c$ is normal and this is
precisely what is asserted in \cite[Corollary 5.6]{Geom II}
\end{proof}

\medskip
For any two points $(s,\lambda), (t, \gamma)$ of $B$ with $s < t$ we
define the {\bf complex slope} of the segment $[(s,\lambda), (t,
\gamma)]$ to be $(\gamma - \lambda)/(t-s)$.  This concept will be
particularly useful if the segment $[(s,\lambda), (t, \gamma)]$ is a face
$F$ of $B$, in which case we call $(\gamma - \lambda)/(t-s)$ the complex
slope of $F$.

\begin{remarks}

\noindent
\smallskip

\noindent (1) Suppose that $\gl$ lies at a corner of $W_k(c)$ and $L_i$
are the sectorial tangent lines of support for $W_k(c)$ at this point.
The faces
$F_1$ and $F_2$ described in part (3) of Theorem 3.4  were specifically
constructed  so that the  corresponding line $\pi_k(L_i)$ in $\mathbb
R\times \mathbb C$ intersect the relative interiors of the $F_i$'s. If
the sharp face $F$ has dimension one, then it bounds precisely two
faces of dimension two in the boundary of $B$.  In this situation, we
say that the union of the
$F_i$'s is a ``shelf\," of $B$ and $F$ is the ``edge" of the shelf.

\medskip
\noindent (2) Now suppose  $F_1 \cap F_2 = F = \{(k/n,k\gl/n)\}$ is a
point.  In this case $(k/n,k\gl/n)$ is an extreme point of $B$ and so
it cannot be an interior point of any face of $B$.  Thus, this point is
an ``end point" of any face that contains it. We say that
$(k/n,k\gl/n)$ is a ``mountain peak" in this case . In contrast to the
one dimensional case, there are three possibilities for the local
geometry of $B$ near $(k/n,k\gl/n)$, which are as follows.

{\smallskip \narrower
\noindent (a) Each face of dimension one that ends at $(k/n,k\gl/n)$
is the edge of a shelf in $B$.

\noindent(b) The point $(k/n,k\gl/n)$ bounds at least one shelf in $B$
and at least one  face of dimension one that is not the edge of a
shelf.

\noindent(c) The point $(k/n,k\gl/n)$ does not bound a shelf.

\smallskip}
\noindent

\medskip
\noindent  If $N$ is abelian so that $B$ has a finite number of
extreme points, then $(k/n,k\gl/n)$ lies at a mountain peak of type
(a).  Observe that in this case, the faces that intersect the lines
$\pi_k(L_i)$ do not intersect in a shelf as in remark (1) above.
However, their boundaries are edges of shelves formed from one of the
$F_i$'s and another two dimensional face of $B$. In most abelian
examples the corners of the isotrace slices all lie in shelves and
there are no mountain peaks.  On the other hand, it is fairly
straightforward to construct abelian examples where this type (a)
phenomenon occurs and it seems possible that there are noncommutative
examples of which also display this geometry.  Examples also show
that mountain peaks of type (b) can occur.

If the boundary of the numerical range is $c$ is nonsingular (i.e.,
does not contain any corners), then $\Psi(0)$ is a mountain peak of
type (c).  Further, it follows from the symmetry of $B$ that $\Psi(1)$
is also a mountain peak of type (c) in this case. We have not been able to
construct an example of  type (c) mountain point at any other points of $B$.

\medskip

\noindent (3) If $F$ is {\em any } sharp face of dimension one in $B$,
then it is the edge of a shelf and if $(k/n,k\gl/n) \in F$ for some
$k$, then $\gl \in W_k(c)$.   To see this, it is convenient to view 
$B$ as a subset of
\rn{3} so that we may use the notation developed in \cite{AAW} and 
\cite{Geom II}. Observe
that  since $F$ is a sharp face of dimension one there must be 
linearly independent
spectral pairs $\spk{i}, i = 1,2$ such that the corresponding faces
\[
F_i = \Psi([\ptspp{i},\ptspp{i}])
\]
have dimension two and contain $F$.  Further, it follows that the
vectors $\tv_i$ are linearly independent.  Indeed, if they were
linearly dependent, then replacing $(s_i,\tv_i)$ with a multiple of
itself if necessary, we could assume that $\tv_1=\tv_2$.  (Recall that
replacing a spectral pair by a multiple of itself leaves the associated
face unchanged).  In this case since $(s_1,\tv_1)$ and $(s_1,\yv_1)$ are
linearly independent, we must have that $s_1 \ne s_2$.  But then, since
these vectors are each normal vectors for $F$ we would get that $F$ has
dimension zero by \cite[Corollary 4.9(2)]{Geom II}, a contradiction.
With this we may repeat the remaining portion of the proof of part (3)
of  Theorem 3.2 verbatim to get that the assertions in part (3) hold
for $F$.

\medskip

\noindent (4)  Now  suppose that  $F = \Psi([z^+,z^-])$ is the sharp
edge of a shelf in $B$ so that $z^\pm$ are central projections in $N$
with $z^- < z^+$.  In this case if we write $z = z^+-z^-$, then we get
that $zb_i =\gb_i z$ for $i =1,2$ and if we  put $\gm = \gb_1 +
i\gb_2$, then $zc = \gm z$ so that
$\gm$ is a reducing eigenvalue of $c$.   This number is the ``complex
slope" of $F$  because
\[
\frac{\tau(cz^+) - \tau(cz^-)}{\tau(z^+) - \tau(z^-)} =
\frac{\tau(zc)}{\tau(z)} =
\frac{\tau(\gamma z)}{\tau(z)} = \gamma.
\]
We note that if $c$ is normal so that every face of dimension one is
the sharp edge of a
shelf in  $B$ and the complex slopes of these faces lie in the
spectrum of $c$, then these
complex slopes fill out the spectrum.  A proof of fact this will
appear as part of a
more general result in a forthcoming paper.

\noindent (5) If  $F = \Psi([z^+,z^-])$ is the sharp edge of a shelf as
in remark (4) and if the rank of $z = z^+-z^-$ is $r$,  then $F$ intersects $r$
adjacent isotrace slices. Thus, in general, the point $(k/n,k\gl/n)$ may be
an interior point of the face $F$.  However, if the reducing eigenvalues of
$c$ have multiplicity one, then $F$ lies between two adjacent isotrace
slices and
$(k/n,k\gl/n)$ is one of the end points of $F$ and so is an extreme point
of $B$.

\end{remarks}

\bigskip

\section{Algebraic and Analytic Preliminaries }

In this section we review the algebraic and analytic facts which can be
used to establish the analyticity of the eigenvalues of matrix pencils.
(See \S 5.3 below for the precise statement of this fact).  Since the
material discussed here is classical the presentation is  brief.
Readers may wish to skim this section to grasp the notation introduced
here and then read Theorem 4.5,  upon which the next section is based.

Suppose $f$ is a complex  polynomial of degree $n$ of the form
\[
f(z,y) = p_0(z)+p_1(z)y +\cdots +(-1)^ny^n.
\]
Since the ring of polynomials $\mathbb C[z,y]$ is a unique
factorization domain \cite[Page 127, Corollary 2]{Hers} and the
coefficient of $y^n$ is a nonzero constant there are irreducible
polynomials $f_1,\dots f_p$ of positive degrees
$d_1,\dots,d_p$ and positive integers $n_1,\dots,n_p$ such that
$d_1n_1+\cdots + d_kn_k = n$ and
\[
f(z,y) = f_1(z,y)^{n_1}\cdots f_p(z,y)^{n_p}.
\]
We call this factorization the {\bf irreducible decomposition} of $f$ in
$\mathbb C [z,y]$.

Recall that if $p(z)$ is a polynomial in a single variable of degree
$n$ and $\gl_1,\dots,\gl_n$ are its roots, then the {\em discriminant}
of $p$ is by definition
\[
\dis(p) = \Pi_{i\ne j}(\gl_i-\gl_j).
\]
Thus, $p$ has a repeated root if and only if $\dis(p) = 0$

Given a polynomial $f(z,y)$ in two variables, we adopt the convention
that the roots of $f(z,y)$ are the roots of the polynomial $p_z$
determined by the formula $p_z(y) = f(z,y)$.  The discriminant of
$f$ is the function of $z$  defined by
\[
\dis(f)(z) = \dis(p_z).
\]
The basic properties of $\dis(f)$ are as follows.

\begin{theorem}The following statements hold.
\begin{enumerate}

\item $\dis(f)$ is a polynomial in $y$.

\item The polynomial $f$ has a repeated factor in its irreducible
decomposition if and only if $\dis(f)$ is identically zero.

\end{enumerate}
\end{theorem}
\begin{proof} This is well known. See \cite[Chapter V
\S10, Corollary, Page 138 and Proposition 5, page 139]{Lang}, for
example.
\end{proof}

\begin{corollary} If $f(z,y)$ is a polynomial of degree $n$ such that
$f(z,y) = f_1(z,y)\cdots f_k(z,y)$, where the
$f_i$'s are irreducible, then there is a finite set $C_f$ of complex
numbers such that if  $z_0\not \in C_f$ where $z_0$ is held constant, then
$f(z_0,y)$ has
$n$ distinct roots.
\end{corollary}
\begin{proof}Since $f$ has no repeated factors, its discriminant is
not identically zero.  Since the discriminant is a polynomial in $y$,
it has a finite number of roots.  Let  $C_f$ denote the roots of
$\dis_f$.  If $z_0\not\in C_f$, then $\dis(f)(z_0)\ne 0$ and so $f(z_),y)$
has $n$ distinct roots.
\end{proof}
The roots of $\dis(f)$ are called the {\bf critical points} of $f$.
More generally, if $f(z,y) = f_1(z,y)^{n_1}\dots f_p(z,y)^{n_p}$ is the
irreducible decomposition of $f$, then we define the critical points of $f$
to be the set of critical points of the {\bf reduced polynomial }
$f_1(z,y)\cdots f_p(z,y)$.


Let $f(z,y)$ denote a complex irreducible polynomial of degree $n$ as
discussed above. An analytic function $\gl$ which  satisfies $f$ in
the sense that
$f(z,\gl(z)) = 0$ for each $z$ in the domain of $\gl$ is called an {\bf
algebraic analytic function}.  Such functions were thoroughly analyzed
by the first half of the last century.  The basic facts of the theory
are summarized below.

\begin{theorem} If  $f$ is an irreducible polynomial of degree $n$ and
$C_f$ denotes its (finite) set of critical points, then  the following
statements hold.
\begin{enumerate}

\item The equation $f(z,y) = 0$ defines precisely one $n$--valued
function $y= F(z)$  which is analytic on $\mathbb C_f = \mathbb C
\setminus C_f$. \item If $D$ is a closed disk in $\mathbb
C_f$ then the $n$ function elements of $F$ determine algebraic
analytic functions  $\gl_1,\dots,\gl_n$ on $D$ such that
\[
f(z,\gl_i(z)) = 0,\quad i = 1,\dots, n.
\]

\item If $D_1,\dots,D_k$ are closed disks in $\mathbb
C_f$ such that $D_i\cap D_{i+1}\ne \emptyset, i = 1,\dots k-1$, and
$\gl_1,\dots,\gl_n$ are the functions  on $D_1$ described in $(2)$, then
each $\gl_i$ has a  unique   analytic  continuation from $D_i$ to $D_{i+1}$

\item If  $D_k = D_1$ and $D_1\cup D_2\cup \cdots \cup D_k$  encircles a
critical point, then  when we return to $D_1$, the same functions are
obtained except that they may have been permuted.
\end{enumerate}
\end{theorem}
See \cite[vol II, chapter 5]{Knop} or \cite[\S 12.3]{San} for a complete
discussion of algebraic analytic functions.   More specifically, part
$(1)$ of  Theorem 4.1 is stated in \cite[vol. II, page 121]{Knop} and
the remaining parts  are proved in \cite[vol. II, \S 15]{Knop}.

       Now suppose  $P$ is a piecewise linear path consisting of lines
joining the
critical points and a half-line  such as the one shown below.

\vskip .2in
\setlength{\unitlength}{1 in}
\begin{picture}(2,2)(-2.5,0)
\put(-.97,-.03){\line(1,2){.225}}
\put(-.55,.03){\line(-1,2){.195}}
\put(-.55,.03){\line(1,2){.34}}
\put(.245,-.2){\line(-1,2){.455}}
\put(.25,-.2){\line(1,2){.35}}
\put(.9,-.1){\line(-1,2){.3}}
\put(.9,-.1){\line(1,2){.407}}
\put(1.65,.03){\line(-1,2){.342}}
\put(-.97,-.03){\vector(-1,1){1}}
\put(-1.8,1){(to $\infty$)}
\put(.21,-.25){$\bullet$}
\put(.86,-.15){$\bullet$}
\put(-1.01,-.076){$\bullet$}
\put(-.59,-.02){$\bullet$}
\put(1.6,0){$\bullet$}
\put(-1.03,-.2){$z_1$}
\put(-.6,-.13){$z_2$}
\put(.2,-.35){$z_3$}
\put(1.6,-.11){$z_m$}
\put(0,1){$P$}
\end{picture}
\vskip .8in

Since $\mathbb C \setminus P$ is simply connected  there are $n$
distinct analytic functions $\gl_1,\dots, \gl_n$ on $\mathbb C
\setminus P$ such that $f(z,\gl_i(z)) = 0$ by parts $(2)$ and $(3)$ of
Theorem 5.3 and  the Monodromy Theorem.

These results immediately yield yield the following facts about the
general case.

\begin{theorem}If $f$ is a polynomial of degree $n$  such that
\[
f(z,y) = f_1(z,y)^{n_1}\cdots f_p(z,y)^{n_p},
\]
where each $f_i$ is irreducible with degree $d_i$, $C_f$ denotes
the critical points of  $f$ and $P$ is any path as
described above, then there are distinct functions
\[
\gl_{lm},\quad 1\le l\le  d_m,1 \le  m\le p
\]
such that:
\begin{enumerate}

\item Each $\gl_{lm}$ is analytic on $\mathbb C\setminus P$.

\item We have $f(z,\gl_{lm}(z)) = 0$ for each $l,m$ and $z \in \mathbb
C\setminus P$.

\item If $z\in \mathbb C\setminus P$ and $(r,s) \ne (l,m)$, then
$\gl_{rs}(z) \ne \gl_{lm}(z)$.
\end{enumerate}
\end{theorem}
We say that functions satisfying property $(3)$ in the Theorem above
are {\bf completely distinct}.


Recall that the {\bf pencil} of matrices determined by $b_1$ and $b_2$
is the set of all matrices of the form $b_z = b_1 + zb_2$, for $z\in
\mathbb C$.  The theory of matrix pencils has a long and honorable
history which includes fundamental contributions by  Kronecker among
others. (See \cite[Vol. II, Chapter XII]{Gant}, for example). In order
to analyze the eigenvalues of $b_z$, write
\[
f(z,w) = \det(b_z-y1) = \det(b_1 + zb_2 -y1).
\]
Since $b_1$ and $b_2$ are $n\times n$ matrices, $f$ is a polynomial of
total degree $n$. Now write write
\[
f(z,y) = f_1(z,y)^{n_1}\cdots f_p(z,y)^{n_p},
\]
for irreducible decomposition of $f$. If we denote the degree of each
$f_i$ by $d_i$ and set $d = d_1+\cdots + d_p$, then we may apply
Theorem 4.4 to get the following result.

\begin{theorem} If $C_f$ denotes the set of critical points of $f$,  $P$
is a path in $\mathbb C$ as described above, and $\mathbb C_P = \mathbb
C\setminus P$, then  there are $d$ completely distinct functions
\[
\gl_{lm},\quad 1 \le l \le d_m,\quad 1 \le m \le p
\]
which are analytic on $\mathbb C\setminus P$ and satisfy $f$.  Thus,
the values of these functions  at a point $z\in \mathbb C\setminus P$
are the $d$ distinct eigenvalues of $b_z$.  Each of the eigenvalues $
\{\gl_{lm}(z): 1\le l\le m\}$ each has multiplicity $n_m$.
\end{theorem}

\section{The Structure of the Boundaries of the Isotraces}

As noted in part 5 of Theorem 0.1, it was shown by Agler and
Narcowich in \cite{Ag} and \cite{Nar}  that if $c$ is a compact
operator and 0 is in the interior of $W(c)$, then the boundary of $W(c)$
is the union of a finite number of analytic arcs. Our goal
in this section is to derive analogous results for the boundary of the
$k$-numerical range.  We note that the numerical range was also studied
by Kippenhahn in \cite{Kip},  where he showed that the boundary of
$W(c)$ is an algebraic curve \cite[Satz 10, page 199]{Kip}.

The basic idea in the proof of these results, which is the same in all
three papers,  is as follows.  If we first translate so that 0 is in the
interior of  $W(c)$, then each ray emanating from the origin determines
a unique tangent line to the boundary of $W(c)$; namely the tangent
line which is perpendicular to the given ray.  It turns out that if
$L_\ta$ denotes the tangent line determined by the ray making the angle
$\ta$ with the
$x$--axis, then its distance $r_\ta$ from the origin is precisely the
maximum eigenvalue of $b_\ta = \cos\ta b_1 + \sin\ta b_2$.  This latter
quantity was then analyzed with the characteristic polynomial in the
finite dimensional case or by a generalization in the infinite
dimensional case.  We will show below that this procedure  easily
generalizes  to include the $k$--numerical range.  Let us
now discuss this in more detail.

Let $C$ denote a compact convex subset of \rn{2} that contains 0 in its
interior.  If we generalize the notation developed above to apply in 
this new setting,
then we have the  following picture.

\vskip .2in
\setlength{\unitlength}{1 in}
\begin{picture}(2,2)(-1.8,0)
\qbezier(-.5,.5)(1,1)(2,0)
\put(-.5,0){\line(1,0){2}}
\put(0,-.5){\line(0,1){2}}
\put(0,0){\vector(1,2){1}}
\put(.2,1.049){\line(2,-1){2}}
\put(.15,.07){$\ta$}
\put(.42,.87){$\bullet$}
\put(.6,.9){$P_\ta = r_\ta e^{i\ta}$}
\put(.7,.3){$C$}
\put(2.25,0){$L_\ta$}
\end{picture}
\vskip .8in

Thus, $L_\ta$ is the tangent line to the boundary of $C$ which is
perpendicular to the ray making the angle $\ta$ with the $x$--axis,
$P_\ta$ is the intersection point of these lines and $r_\ta$ is the
distance from $P_\ta$ to the origin.

\begin{lemma}[Agler] The function $\ta \to r_\ta$ is differentiable if
and only if
\[
\lim_{\phi\to\ta}L_\phi\cap L_\ta = P_\ta.
\]
Further, the points of non--differentiability of $r_\ta$ are in
one-to-one correspondence with the line segments in the boundary of $C$.
\end{lemma}
\begin{proof}This is basically due to Agler \cite{Ag}, although it is
not explicitly stated there.   He observed that if $\phi\ne \ta$,
then
\begin{equation}
L_\ta\cap L_\phi = \left(\frac{r_\phi\sin\ta
-\gat\sin\phi}{\sin(\ta-\phi)},\frac{r_\phi\cos\ta
-\gat\cos\phi}{\sin(\ta-\phi)}
\right). \tag{$*$}
\end{equation}
(The formula presented in \cite[(4.2)]{Ag} appears to
contain a typographical error since it differs from the one offered
above by a minus sign in the second coordinate).  Hence, if
$\lim_{\phi\to\ta}L_\phi\cap L_\ta = P_\ta$, then $r_\ta$ is
differentiable, and if $r_\ta$ is differentiable, then it follows from
$(*)$ that this limit exists.

Next, Agler showed   in \cite[Lemma 4.2]{Ag} that
if
\[
L_\ta\cap C = \{tP_1+(1-t)P_2:0\le t\le 1\}
\]
where $P_1$ comes from moving counterclockwise around $\partial C$, then
\[
P_1 = \lim_{\phi\to \ta^-}L_\phi\cap L_\ta  \text{ and }
P_2 = \lim_{\phi\to \ta^+}L_\phi\cap L_\ta.
\]
Hence, if the boundary of  $C$ contains a line segment with end points
$P_1\ne P_2$, then $\lim_{\phi\to \ta} L_\phi\cap L_\ta$ does not exist
and so $r_\ta$ is not differentiable by the first part of the proof.
Conversely, if $r_\ta$ is not differentiable, then  $\lim_{\phi\to \ta}
L_\phi\cap L_\ta$ does not exist and so
\[
P_1 =  \lim_{\phi\to \ta^-}L_\phi\cap L_\ta
       \ne \lim_{\phi\to \ta^+}L_\phi\cap L_\ta = P_2.
\]
Since $P_1\ne P_2$, the line segment joining $P_1$ and $P_2$ lies on
the boundary of $C$.
\end{proof}

\medskip

Now suppose that  $C = W_k(c)$ and $0$ is in the interior of $W_k(c)$.
By part$(5)$ of Theorem 0.3, we have that if $\gb_k^+$ is the sum of the
$k$ largest eigenvalues of $\text{Re}(c)$,  then the line $x=\gb_k^+$ is
tangent to $W_k(c)$.   If we fix $0\le \ta < 2\pi$, then applying
this argument to $W_k(e^{-i\ta}c)$ we get that the line $x = \gb_{\ta,k}^+$
is tangent to $W_k(e^{-i\ta}c)$, where $\gb_{\ta,k}^+$ is the sum of the
$k$ largest (not necessarily distinct) eigenvalues  of
\[
b_\ta = \text{Re}(e^{-i\ta}c) = \cos\ta b_1 + \sin\ta b_2.
\]

If we now rotate back by the angle $\ta$, we get that the line $L_\ta$
is at distance $\gb_{\ta,k}^+$ from the origin. That is,
\[
L_\ta = \{\gb_{\ta,k}^+e^{i\ta} + te^{i(\ta +\pi/2)}:t\in \mathbb R\}
\]
is tangent to $W_k(c)$.  In other words, $r_\ta = \gb_{\ta,k}^+$ in this
case.  We now record this fact for future reference.

\begin{theorem}If $\gb_{\ta,k}^+$ denotes the sum of the $k$ largest
(not necessarily distinct) eigenvalues of $b_\ta$, then
$r_\ta = \gb_{\ta,k}^+$.
\end{theorem}

Now write
\[
g(u,v,w) = \det(ub_1 + vb_2 -w1)
\]
so that for each fixed pair $(u,v)$, $g(u,v,w)$ is the characteristic
polynomial of $ub_1 + vb_2$.  Thus, the eigenvalues of $b_\ta$ are
precisely the roots of
\[
       g(\cos\ta,\sin\ta,w) = \det(\cos\ta b_1 + \sin\ta b_2 -w1).
\]

As $g$ is a homogeneous polynomial, we may write $\displaystyle z =
\frac{v}{u}\text{ and } y = \frac{w}{u}$ and set $f(z,y) = g(1,z,y)$.
Note that  we have the relation
\begin{equation}
u^nf(v/u,w/u) = u^ng(1,v/u,w/u) = g(u,v,w). \tag{$**$}
\end{equation}
Also, we have
\[
f(z,y) = \det(b_1+zb_2-y1).
\]
and so we may apply the results described in Theorem 4.5 to $f$.

We now translate these results to obtain analogous results for the
eigenvalues of $b_\ta$ as follows.  If $\gb_\ta$ is an eigenvalue for
$b_\ta$, then we have $g(\cos\ta,\sin\ta,\gb_\ta) = 0$ and so using
the relation $(**)$ above, we get that if $\cos\ta \ne 0$, then
\[
0= g(\cos\ta,\sin\ta,\gb_\ta) = \cos^n\ta g((1,\tan\ta,\sec\ta\gb_\ta)=
\cos^n\ta f(\tan\ta,\sec\ta\gb_\ta).
\]
Hence, applying Theorem 4.5, we get
\[
\sec\ta \gb_\ta = \gl_{ij}(\tan\ta),\quad \gb_\ta =
\cos\ta\gl_{ij}(\tan\ta)
\]
for some pair $(i,j)$.  Thus, we may write
\[
\phi_{ij}(\ta) = \cos\ta\gl_{ij}(\tan\ta),\quad 1\le i\le d_j,\quad
1\le j\le p
\]
for the eigenvalues of $b_\ta$.

    Further, if we select the path $P$ in part $(3)$ of
Theorem 4.5 so that $P\cap \mathbb R = C_f\cap \mathbb R$ as shown below
\footnote{Since the matrices $b_1$ and $b_2$ are
self--adjoint, the coefficients of the polynomial $\dis_f$ are real and
so the critical points occur in conjugate pairs}

\setlength{\unitlength}{1 in}
\begin{picture}(2,2)(-2.5,0)
\put(-1.5,0){\line(1,0){3.25}}
\put(-.97,-.01){\line(1,2){.213}}
\put(-.55,.0){\line(-1,2){.208}}
\put(-.55,.0){\line(5,-2){.75}}
\put(.2,-.3){\line(0,1){.6}}
\put(.836,0){\line(-2,1){.6}}
\put(.83,-.02){\line(1,2){.425}}
\put(1.253,-.75){\line(0,1){1.6}}
\put(-.95,-.02){\vector(-1,1){1}}
\put(-1.8,1){(to $\infty$)}
\put(.16,.28){$\bullet$}
\put(.16,-.343){$\bullet$}
\put(.805,-.04){$\bullet$}
\put(-1.01,-.04){$\bullet$}
\put(-.59,-.04){$\bullet$}
\put(1.21,.78){$\bullet$}
\put(1.213,-.8){$\bullet$}
\put(-1.03,-.15){$z_1$}
\put(-.6,-.15){$z_2$}
\put(.15,-.44){$z_3$}
\put(.15,.45){$z_4$}
\put(1.21,-.9){$z_m$}
\put(1.8,-.06){$\mathbb R$}
\put(0,1){$P$}
\end{picture}
\vskip .8in
\bigskip

\noindent then we get that there is a finite set $\Theta_f$
in $[0,2\pi)$, which we may assume contains 0,  such that if $\ta_l$ and
$\ta_{l+1}$ denote successive elements in $\Theta_f$ and we denote the open
interval  $(\ta_l,\ta_{l+1}) $ by  $I_l$, then each $\phi_{ij}$ is real
analytic on $I_l$ and the functions $\phi_{ij}$ are distinct on $I_l$.
The results of this analysis are recorded below.

\begin{theorem} There are  positive integers $d_1,\dots, d_p$ and
$n_1,\dots,n_p$ and a a finite subset $\Theta_f$ of $[0,2\pi)$ such the
following statements hold.
\begin{enumerate}

\item We have $d_1n_1+\cdots+d_pn_p = n$.

\item If $\ta \not \in \Theta_f$, then $b_\ta$ has precisely $d$
distinct eigenvalues, where $d = d_1+\cdots +d_p$.  If $\ta \in
\Theta_f$, then $b_\ta$ has fewer than $d$ distinct eigenvalues.

\item For $\ta \not \in \Theta_f$ and $1\le i\le p$, $b_\ta$ has $d_i$
distinct eigenvalues with multiplicity $n_i$.

\item There are completely distinct functions $\phi_{ij}$ for $1\le
i\le d_i$ and $1\le j\le p$ which are real analytic on
$[0,2\pi)\setminus\Theta_f$  and such that $\phi_{ij}(\ta)$ is an
eigenvalue of multiplicity $n_j$ for $b_\ta$.

\end{enumerate}

\end{theorem}

\begin{theorem} If $0 <k < n$, then the boundary of the $k$--numerical
range consists of a finite number of line segments and curved real
analytic arcs.
\end{theorem}
\begin{proof} We use the notation introduced in Theorem 4.3.  Let $I$
denote an open interval in $[0,2\pi)\setminus \Theta_f$ whose end points
are adjacent points in $\Theta_f$. We have then that the functions
$\phi_{ij}$  are  real analytic and completely distinct on $I$ by part
$(4)$ of Theorem 5.3. Hence, we may re-index these functions as $\phi_i$
for $i = 1,\dots, n$ and assume that
\[
\phi_1(\ta) > \cdots > \phi_d(\ta)
\]
for each $\ta \in I$.  Similarly, with additional re-indexing, we may
assume that for each fixed $i$ the eigenvalue $\phi_i(\ta)$ has
multiplicity $n_i$.  Hence, we may  find positive integers $l$ and
$m$ such that $1 \le m \le n_{l+1}$ and
\[
k = n_1 +\cdots +n_l + m.
\]
With this we get that $\gb_{\ta,k}^+$, which is the sum of the $k$ largest
eigenvalues of $b_\ta$ by Theorem 5.2, is precisely
\[
r_\ta =n_1\phi_1(\ta) + \cdots + n_l\phi_l(\ta) + m\phi_{l+1}(\ta)
\]
and so the map $\ta\to \gb_{\ta,k}^+$ is real analytic on $I$ by Theorem
4.3.

Since the set $\Theta_f$ is finite,  this argument shows that
$\gb_{\ta,k}^+$ is differentiable on $[0,2\pi)$, with the possible
exception  of  a  finite subset of $\Theta$.   Hence,the boundary of
$W_k(c)$ contains a finite number of line segments by Lemma 5.1 and
the remaining boundary points are contained in  a finite number of
curved real analytic arcs.

\end{proof}

\begin{corollary} The spectral scale $B(c)$ has a finite number of faces
of dimension two.
\end{corollary}
\begin{proof} If $F$ is a face in $B$ of dimension two, then  $F$ is
transverse to the isotrace slices of $B$ by \cite[Lemma 6.4(1) ]{Geom
II}.  Hence, it must  intersect at least one isotrace in a line segment by
\cite[Lemma 6.4(1)]{Geom II}. Since each isotrace slice contains only a
finite number of line segments by Theorems 1.3 and 5.4, $B$ can have only
a finite number of faces of dimension two.

\end{proof}

\remarks
As noted at the beginning of this section,  Agler, Kippenhahn and
Narcowich each base their analyses on the fact that the tangent line
$L_\ta$ for $W(c)$ is determined by the largest eigenvalue of $b_\ta$.
Kippenhahn then continues his analysis  in the setting of algebraic
geometry. On the other hand,  Agler uses the theory of  Several Complex
Variables in his analysis, while Narcowich bases his proof on a Theorem
of Nagy.

The argument in Agler's paper is quite similar to the one presented
here.  For example, both reduce to a consideration of polynomials
without repeated factors and then  rely on the discriminant to determine
the set $\Theta$.  The main difference between the proof above and
Agler's is that  since we are working in finite dimensions, we
can base our analysis on the  classical results on algebraic analytic
functions, rather than delving into the more abstruse
world (from our vantage, at least) of several complex variables. Our
contribution consists of the observation that Agler's proof may be
extended to cover the case where we are working with sums of
eigenvalues.   Nagy's theorem (\cite[Page 376]{Nagy}, which is the
basis for Narcowich's argument may be viewed in this context as an
infinite dimensional version of  Theorem 4.3.

Kippenhahn showed that the boundary of the numerical range lies in
the zero set of an algebraic curve as follows.  He considered $f(u,v,w) =
\det(ub_1 + vb_2 +w1)$, which is a homogeneous polynomial and so defines
an algebraic curve in complex projective 3-space. He observed if $r_\ta$
denotes the maximum eigenvalue of $\cos\ta b_1 + \sin\ta b_2$ for each
fixed $\ta$, then  the formula $b_\ta =\cos\ta u +\sin\ta v + r_\ta w = 0$
determines the tangent line $L_\ta$  described at the beginning of this
section.  The envelope of these tangents is also  an algebraic curve,
which is called the {\bf dual curve} to $f$ \cite[page 86 and 
Proposition 5, Page
253]{Bris}.  Thus, the boundary of the numerical ranges lies in the 
zero set of this dual
curve.

It seems likely that the boundary of the $k$-numerical range
also lies on an algebraic curve.  To show this using Kippenhahn's
technique requires showing that if $r_\ta$ denotes the sum of the $k$
largest eigenvalues of $b_\ta$, then points of the form
$(\cos\ta,\sin\ta,r_\ta)$ also lie on an algebraic curve.  Now it is a
standard algebraic that fact that since each eigenvalue of $b_\ta$ lies
on a fixed  algebraic curve.  This means that for each fixed $\ta$ we may
find real numbers $u_\ta$ and $v_\ta$ such that point$(u_\ta,v_\ta,
r_\ta)$ lies on this curve. But it seems possible that we would have
$u_\ta,v_\ta) \ne (\cos\ta,\sin\ta)$ for most $\ta$'s.

\endremarks

\section{Conjectures}
\bigskip

Conjecture 6.1. If $1 \le k < n$ and $\lambda$ is a corner of $W_k(c)$,
then
$\phi(\lambda)$ lies on a shelf in B.
\bigskip

\noindent
Conjecture  6.2. If $n$ is even and $I_{n/2}$ is a polygon, $c$ is
normal.
\bigskip

\noindent
Conjecture  6.3.  The boundary of $W_k(c)$ is contained in the real 
zero set of an
algebraic curve.

\end{document}